\documentclass[a4paper,10pt]{article}

\def \Z {{\mathbf Z}}
\def \R {{\mathbf {R}}}
\def \N {{\mathbf {N}}}

\def \B {{\cal B}}

\def\uu{\bigsqcup}
\def\eps{\varepsilon}
\title{ Некоторые задачи о сохраняющих меру преобразованиях}
\author{В.В. Рыжиков}
\date{}
\usepackage[T1]{fontenc} 
\usepackage[cp1251]{inputenc} 
\usepackage[russian]{babel}
\usepackage{graphicx}
\graphicspath{{pictures/}}
\DeclareGraphicsExtensions{.pdf,.png,.jpg}
\begin{document}
\Large

\maketitle
{  \large Предлагается версия рассказа ''Задачи, которые нас выбирают''. Обсуждается факторизация 
преобразований в произведение инволюций, теорема Фюрстенберга о слабом кратном перемешивании в среднем вдоль прогрессий,
теорема Рота-Фюрстенберга о кратном возвращении, кратное перемешивание систем с эргодической группой Гордина, отсутствие свойства Лерера-Вейса для действия решеток. Упомянуты необычные свойства сохраняющих меру конструкций, 
парадоксальный эффект Катка,  некоторые нерешенные задачи и известные проблемы в эргодической теории.
 }

\section{Три инволюции и немного воспоминаний.}
Когда-то не было интернета.  Студенты ходили в библиотеки, посещали читальные залы, заглядывали в книжные магазины.  
Я купил книжку Улама ''Нерешенные математические задачи'' \cite{Ulam}.  
В ней среди прочих была задача: \it является ли группа всех обратимых
сохраняющих меру преобразований отрезка $[0,1]$  алгебраически простой?  \rm Задача  показалась интересной, а ответ
очевидным. Понятно, что сохраняющее меру преобразовние  похоже на перестановку бесконечно малых отрезков.
Перестановку можно считать четной, так как   нечетные отличаются от четных лишь на одну бесконечно малую транспозицию.
 А группа четных перестановок проста. 

Сообщил эту  идею научному руководителю Альберту Григорьевичу Драгалину. 
Он  предложил  мне заняться нестандартным анализом, избавив  от изучения книги 
Клини и Весли ''Основания интуиционистской математики''. 
Хотя нестандартный анализ   не помог, мне
удалось  показать, что любое преобразование является композицией не более, чем 8 инволюций 
(сохраняющих меру преобразований периода 2). Это утверждение -- наиболее важный момент   
в доказательстве простоты группы, остальное -- дело относительно несложной техники.

Возник  естественный вопрос: была ли эта задача из книги Улама решена кем-то ранее?
Преподаватель  фунционального анализа  Евгений Михайлович Никишин советует обратиться к  Анатолию Михайловичу Степину,  знатоку групп преобразований.     На  спецсеминаре узнал,  что скорее всего задача решена, но надо было проверить. На следующий день  в читальном зале   после многочасовой обработки множества увесистых томов  реферативного журнала   нахожу искомый реферат. Задачу   решил  Альбер Фати  \cite{Fathi}.
Удивительно, но и  Фати  использовал результат о факторизации преобразования в произведение 10 инволюций. Его доказательство было очень похоже на мое.

В книге Халмоша \cite{H} был упомянут результат японского математика Анзаи о существовании
преобразований, не сопряженных своему обратному. Следовательно,  двух инволюций недостаточно для факторизации. 
 Во время одной из прогулок в  парке  Дружбы у Китайского посольства
  вдруг  возникла  мысль о том, как  из 8-ми  (аналогично из 88-ми) инволюций сделать 3 инволюции. 
 Позднее я подружился с   Жан-Полем Тувено и  устроил ему экскурсию в  
 ''Parc des trois involutions'', где нас  посетили новые интересные идеи.  

\vspace{5mm}
\begin{center}
\includegraphics{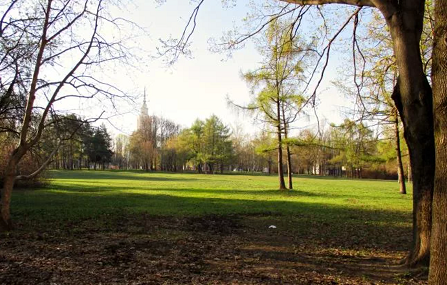}
\newpage
\includegraphics{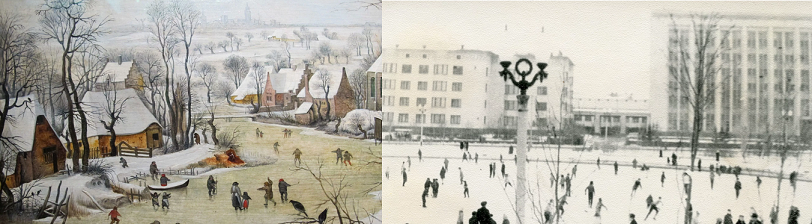}
\end{center}

\begin{center}
\it Справа парк Дружбы в давние времена. \rm
\end{center}

Для доказательства использовалась  лемма  Рохлина-Халмоша в следующей форме:
\it для $T$, обратимого сохраняющего меру апериодического преобразование $[0,1]$,  найдется измеримое множество $B$ такое,
что   $[0,1]$ является объединением  множества
$$B, TB,T^2B,\dots, T^{10}B, \hat B,$$ причем  
$T\hat B\subset B.$ \rm        

       Рассмотрим инволюцию $s$, которая возвращает множество $\hat B$ в $T^{10}B$, а вне 
множеств $\hat B$ и $s\hat B$ является тождественной.
Тогда $sT$ является циклической перестановкой множеств $B, TB,T^2B,\dots, T^{10}B$, которая есть композиция $\tilde TP$,
где $\tilde T$ вне множества $B$ является тождественным преобразованием, а $P$ -- периодическая перестановка множеств
$B, TB,T^2B,\dots, T^{10}B$  ($P^{11}=Id$). 

Пусть мы знаем, что $\tilde T$ является композицией 10 инволюций 
$\tilde S_1, \tilde S_2,\dots, \tilde S_{10}$, которые вне множества $B$ совпадают с  тождественным преобразованием.  Тогда циклическую перестановку $\tilde TP$ можно представить
в виде $S_1S_2\dots S_{10}P$, где инволюции $S_{i+1}$ вне множества $T^iB$ совпадают с  тождественным преобразованием. (Подумайте, как растащить композицию $S=\tilde S_1\tilde  S_2 \dots \tilde S_{10}$ по циклу, где $S_i$ суть копии $\tilde S_i$.)
Получаем факторизацию в три инволюции
$$T= (s \tilde S_1\tilde  S_2 \dots \tilde S_{10})P=SS'S'',$$
где $P=S'S''$ -- почти очевидное представление периодического преобразования как композиции двух ''симметрий'',
$S$  как произведение коммутирующих инволюций является  инволюцией.  

Факторизация преобразования в 10 инволюций требует  больше усилий, так как 
 приходится иметь дело с бесконечными произведениями, нужен трюк, чтобы распутать бесконечно запутанный клубок.
Представьте себе убывающую последовательность вложенных друг в друга множеств $A_i$ (мера $A_{i+1}$ почти в два раза меньше
меры $A_i$), на каждом из них действует инволюция $S_i$ (вне множества  $A_i$ преобразование  $S_i$ является тождественным).
Тогда имеет смысл бесконечное произведение $\prod_i S_i$. После недолгих размышлений становится ясно, что любое преобразование можно  представить в виде такого произведения. Не будем раскрывать секрета, вдруг читатель захочет сам
''победить дракона''. 
Заметка \cite{1985} содержит некоторые  усиления (''всякое преобразование есть коммутатор и композиция двух преобразований периода три''). Результат об инволюциях    был перенесен автором на  булевы алгебры для защиты диплома в 1985 году. После отъезда А.Г. Драгалина в Венгрию  моим руководителем стал Владимир  Андреевич Успенский. 
Модифицированная теорема   звучала так: \it автоморфизм полной булевой алгебры 
является композицией трех инволюций. \rm Есть соответствующая  публикация  с картинками \cite{3}.

\vspace{3mm}
\bf Неразрешимость элементарной теории группы $\bf Aut(X,\B,\mu)$. \rm Выпускник кафедры математической логики Александр Митин доказал, что элементарная теория группы обратимых сохраняющих меру преобразований отрезка $[0,1]$
неразрешима \cite{Mitin}, погрузив  в нее арифметику. Сказанное  означает, что нет алгоритма, который сообщал  бы нам, истинны или ложны для нашей группы высказывания  типа ''существует преобразование, не 
сопряженное  своему обратному''.  Недавно появилась   работа М. Формана \cite{Fo}, в которой автор интерпретирует разнообразные  математические проблемы  в виде  задач об изоморфизме преобразования своему обратному. В маленьком блюдце
отражается вселенная.

\Large


\section{Самоприсоединения и кратная возвращаемость.}
В 1985 году, читая книгу Грэхема ''Теория Рамсея'', с  удивлением узнал, что знаменитая теорема Семереди  (он, кстати,  учился на мехмате  МГУ) о прогрессиях эквивалентна тому, что  для любого множества $A$  положительной меры и обратимого сохраняющего меру преобразования $T$ вероятностного пространства найдется $i>0$, 
для которого  выполняется свойство кратного возвращения:
$$  \mu( A\cap T^iA\cap T^{2i}A\dots \cap T^{ki}A)\ >0.$$
Фюрстенберг нашел  динамическое    доказательство теоремы Семереди \cite{F}.
Случай $k=1$ давно известен как  теорема Пуанкаре о возвращении, которая очевидна, но замечательна. 
 При $k>1$ и особенно при $k>2$ доказательство становится чрезвычайно нетривиальным. 

\vspace{3mm}
\bf Кратное слабое перемешивание в среднем вдоль прогрессий. \rm

Чтобы установить упомянутый  результат  методами эргодической теории, 
Фюрстенберг использовал вспомогательную теорему: \it 
для слабо перемешивающего преобразования $T$ вероятностного 
пространства   $(X,\B,\mu)$ для любых $ A,A_1,\dots A_k\in\B$ при $N\to\infty$ выполнено
$$  \frac 1 {N}  \sum_{i=1}^{N} \mu( A\cap T^iA_1\cap T^{2i}A_2\dots \cap T^{ki}A_k)\ \to 
\ \mu(A)\mu(A_1)\dots\mu(A_k). \eqno {\bf (a.p. Mix(k))}$$\rm

Я попробовал   доказать  это утвержение, не подглядывая в первоисточники. Мной был придуман следующий метод. 
Рассмотрим случай $k=2$. Из любого бесконечного подмножества натурального ряда 
всегда можно выбрать подпоследовательность $N(k)\to\infty$ так, чтобы для всех $ A,A_1,\dots A_k\in\B$ 
$$  \frac 1 {N(k)}  \sum_{i=1}^{N(k)} \mu( A\cap T^iA_1\cap T^{2i}A_2)\ \to \ \nu(A\times A_1\times A_2),
$$
где  $\nu(A\times A_1\times A_2)$ -- это пока такое своеобразное обозначение предела выражений слева. Несложно заметить, что 
 $\nu$ как функция  на цилиндрах $A\times A_1\times A_2$, образующих полукольцо,  является мерой.
Очевидна ее инвариантость 
$$\nu(A\times A_1\times A_2) =\nu(TA\times TA_1\times TA_2),$$
а в силу усреденения имеет место дополнительная инвариантность 
$$\nu(A\times A_1\times A_2) =\nu(A\times TA_1\times T^2A_2).$$
Мера $\nu$ как мера на кубе проектируется на сомножители в меру $\mu$.  Такие $T\times T\times T$-инвариантные меры с хорошими проекциями    называется самоприсоединениями ( self-joinings) или джойнингами.  
Ревнители чистоты русского языка без энтузиазма  воспринимают   фонетическую кальку  ''джойнинг'', поэтому лучше использовать слово  ''присоединение''. 

Свойство слабого перемешивания эквивалентно  эргодичности преобразования $T\times T^2$, которое в свою очередь означает
равносоставленность  можеств $ A_1\times A_2$ и  $ B_1\times B_2$ при условии $ \mu(A_1)\mu(A_2)=\mu(B_1)\mu(B_2).$ 
Последнее означает, что  
$$ A_1\times A_2\ =\uu_{n=1}^\infty (C_n\times D_n), $$
$$ B_1\times B_2\ =\uu_{n=1}^\infty (T^{p(n)}C_n\times T^{2p(n)}D_n),$$
где равенства выполняются с точностью до множеств нулевой $\mu\times\mu$-меры.
 Получаем
$$\sum_n \nu(A\times C_n\times D_n)=\sum_n \nu(A\times T^{p(n)}C_n\times T^{p(n)}D_n),$$
$$\nu(A\times A_1\times A_2) =\nu(A\times B_1\times B_2).$$
Из сказанного вытекает, что  $\nu=\mu\times\mu\times\mu$.
А это означает
$$  \frac 1 {N}  \sum_{i=1}^{N} \mu( A\cap T^iA_1\cap T^{2i}A_2)\ \to \ \mu(A)\mu(A_1)\mu(A_2).$$
Далее по индукции можно установить ${\bf (a.p. Mix(k))}$ для всех $k>2$.

Таким образом, автор своим способом  доказал  упомянутое утверждение  Фюрстенберга, 
 попутно переоткрыв  джойнинги. Фюрстенберг  давно  ввел это понятие
в эргодическую теорию, но  не использовал   для своего  доказательства  ${\bf (a.p. Mix(k))}$, применяя  более сильные методы. 

\vspace{4mm}
\bf Эргодический аналог теоремы Рота о прогрессиях. \rm
Сейчас мы расскажем, опуская некоторые  детали,  как  устанавливается возвращение кратности $2$,
следуя версии    Тувено (Thouvenot) и автора (the author).

\vspace{2mm}
\bf  Теорема (Roth, Furstenberg). \it Для любого множества $A$  положительной меры 
и обратимого сохраняющего меру преобразования $T$ вероятностного пространства $(X,\B,\mu)$
найдется $i>0$, для которого  выполняется 
$$  \mu( A\cap T^iA\cap T^{2i}A\ >0.$$\rm

Доказательство. Благодаря теореме Рохлина о разложении инвариантной меры на эргодические 
компоненты  общий случай сводится к случаю эргодического преобразования $T$. Пусть $Tf(x)$ обозначает $f(Tx)$.
Отметим, что эргодичность оператора $T$ эквивалентна выполнению  условия 
$$\frac 1{N}\sum_{i=1}^{N} T^{i}f\to_w const =\int_X f\ d\mu, $$
для всех $f\in L_\infty(\mu)$. 

Определим оператор $J:L_2(\mu)\to L_2(\mu\times\mu)$ формулой
$$(Jf,g\otimes h)=\lim_{N_k}\frac 1{N_k} \sum_{i=1}^{N_k}\int_X f\ T^ig\ T^{2i}h\ d\mu,$$ 
где $f,g,h\in L_\infty(\mu)$.
 Оператор $J$ корректно определен, так как в силу     эргодичности преобразования $T$ мы получим
$$\lim_{N_k}\frac 1{N_k} \sum_{i=1}^{N_k}\int_X  T^ig\ T^{2i}h\ d\mu   =
\lim_{N_k}\frac 1{N_k}\sum_{i=1}^{N_k}\int_X  g\ T^{i}h\ d\mu=\int_X  g\ d\mu \int_X h\ d\mu.$$
    Усреднение, фигурирующее в  определении оператора  $J$,  обеспечивает важное для нас  равенство
$$(T\otimes T^2)J=J.$$

Замечание. Появившееся равенство  непосредственно связано с  обсуждавшейся дополнительной инвариантностью 
$(Id\times T\times T^2)\nu=\nu$  для самоприсоединения $\nu$, отвечающего оператору $J$:
$$\nu(A\times B\times C) =(J\chi_A, \chi_B\otimes \chi_C).$$  Самоприсоединения и сплетающие операторы 
 в ряде случаев суть две стороны одной медали, выбор сторон диктуется соображениями удобства.

Вернемся к доказательству.  Так как   $(T\otimes T^2)Jf=Jf$, образ оператора $J$ состоит из функций,  неподвижных относительно оператора $T\otimes T^2$. Но пространство таких неподвижных функций является подпространством 
$V\otimes V$, где $V$ порожденное всеми собственными векторами оператора $T$. 
Сказанное  вытекает из несложных  стандартных фактов спектральной теории унитарных операторов, 
подробности  опускаем.

Обозначим через $\pi$ ортопроекцию $L_2(\mu)$ на $V$. Тогда
$$(T\otimes T^2)Jf=Jf, \ \ (\pi\otimes\pi)J=J,$$
$$ (Jf,g\otimes h)= (Jf,(\pi\otimes\pi)g\otimes h)=
\lim_{N_k}\frac 1{N_k} \sum_{i=1}^{N_k}\int_X f\ T^i\pi g\ T^{2i}\pi h\ d\mu =$$
$$=\lim_{N_k}\frac 1{N_k} \sum_{i=1}^{N_k}\int_X \pi f\ R^i\pi g\ R^{2i}\pi h\ d\mu,$$
где $R$ -- ограничение $T$ на $V$.  Но степени оператора $R$  регулярно оказываются близкими к  тождественному
оператору, иначе говоря,  множество тех $i$, когда 
$$\|R^i\pi f -\pi f\|<\eps,  \ \|R^{2i}\pi f -\pi f\|<\eps$$
имеет положительную плотность.
Положив $f=g=h=\chi_A$ при условии $\mu(A)>0$, получаем  
$$(J\chi_A,\chi_A\otimes \chi_A)=\lim_{N_k}\frac 1{N_k} \sum_{i=1}^{N_k}
\int_X \pi \chi_A\ R^i\pi \chi_A\ R^{2i}\pi\chi_A\ d\mu \ >\ 0.$$
Тем самым мы показали, что неравенство   $ \mu( A\cap T^iA\cap T^{2i}A\ >0$ выполняется 
для бесконечного множества значений $i$, что завершает доказательство.

\section{Проблема Рохлина. Гармоническая группа и бесконечные  нити. Гомоклиническая группа}
 \rm Самоприсоединения оказались     полезными  для изучения свойства кратного перемешивания. 
В эргодической теории  с 1949 года
открыт следующий вопрос Рохлина: \it если преобразование перемешивает, будет ли оно перемешивать
кратно? \rm  Дадим определения, следуя \cite{Rohlin}.

\it Преобразование $T$ на вероятностном пространстве назвается перемешивающим, если для всяких  измеримых $A,A_1$ 
имеет место сходимость
$$\mu(A \cap T^{i}A_1)\ \to \ \mu(A)\mu(A_1), \ i\to\infty.  \eqno  {\bf Mix(1)}$$

Преобразование $T$ на вероятностном пространстве называется перемешивающим кратности 2, 
если для всех измеримых множеств $A,A_1,A_2$ выполнено
$$\mu(A \cap T^{i}A_1\cap T^{i+j}A_2)\ \to \ \mu(A)\mu(A_1)\mu(A_2), \ \ i,j\to +\infty.\eqno  {\bf Mix(2)}$$ 
\rm

Поясним, как  присоединения  связаны с проблемой Рохлина ${\bf Mix(1)=? Mix(2)}$.
Если есть гипотетический контрпример $T$, то для некоторых последовательностей 
$ i_k,j_k\to +\infty$ получим
$$\mu(A \cap T^{i_k}A_1\cap T^{i_k+j_k}A_2)\ \to \ \nu(A\times A_1\times A_2)\neq \mu(A)\mu(A_1)\mu(A_2).$$
Мера $\nu\neq \mu\times\mu\times\mu$ называется нетривиальным самоприсоединением cо свойством парной независимости
(иначе говоря, проекции джойнинга $\nu$  на двумерные грани куба $X\times X\times X$ суть меры $\mu\times\mu$). Иногда удается доказать, что таких нетривиальных объектов у изучаемой перемешивающей динамической системы нет, что
 влечет за собой свойство кратного перемешивания. 
 
 Узнать историю   результатов о кратном перемешивании, в том числе полученных благодаря исследованию  самоприсоединений и тесно связанных  с ними    марковских сплетающих операторов,  можно  в работе  \cite{diss}, снабженной соответствующим введением.

{\large  \bf Пример Ледрапье.  \rm Для действий группы $\Z^2$ имеются замечательные примеры перемешивающих действий, не обладающих кратным перемешиванием. Они реализуются коммутирующими автоморфизмами некоторых компактных коммутативных групп.
В оригинальном примере Ледрапье фигурировала  группа, состоящая из 
     функций $h: \Z^2\to\Z_2$,  удовлетворяющих условию гармоничности:

$$h(x, y)\ = \ h(x + 1, y) + h(x,y + 1) + h(x - 1, y) + h(x, y - 1). \ \ (mod 2)$$

\bf  Загадочная доминирующая нить. \rm С гармоническими  последовательностями  связан необычный эффект, который автор обнаружил в 80-х годах, когда появились персональные компьютеры.  Если нарисовать случайную гармоническую последовательность (значению 1 соответствует белый цвет), то отчетливо видны нити.  Возникла гипотеза  о том,  что в типичном случае присутствует   бесконечная  нить (компьютер это явно подсказывает), которая захватывает около   половины  плоскости.  
Но как это доказать? Начнешь думать об одной задаче,   наткнешься на  другую, и обе останутся нерешенными. И третья тут как тут.

\begin{center}
\includegraphics{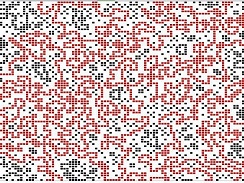}
\end{center}
\large }

\it Встанем на гармоническую нить и пойдем по ней. Возикает бесконечная последовательность из трех символов: 
-1 (налево), +1 (направо), 0 (вперед). Что можно сказать
о свойствах $\Z$-действия, отвечающего такой символической системе?  \rm
\begin{center}
\includegraphics{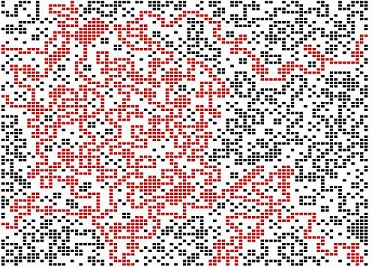}
\end{center}

\vspace{3mm}
\bf Гомоклиническая группа и кратное перемешивание. \rm В качестве  примера применения самоприсоединений для обнаружения 
свойства кратного премешивания приведем  обобщение одного из утверждений работы  \cite{G}.  

Подмножество $H(T)$ группы $Aut(X,\B,\mu)$ всех сохраняющих меру 
преобразований, 
$$ H(T)=\{S\in Aut(X,\B,\mu):  T^{-n}ST^n\to Id,\ n\to\infty \},$$
называется \it гомоклинической группой Гордина. \rm Для $k=1$  в  работе \cite{G} 
 дано доказательство следующего утверждения:

\it Если группа $H(T)$ преобразования $T$ эргодична, то  $T$ обладает свойством ${ \bf Mix(k)}$.\rm

Считая, что ${ \bf Mix(1)}$ установлено, докажем ${ \bf Mix(2)}$. Пусть для некоторых
 $i_m,j_m\to\infty$, $j_m-i_m\to\infty$, для любых множеств $A,A_1,A_2\in\B$  выполнено
$$\mu(A \cap T^{i_m}A_1\cap T^{j_m}A_2)\ \to \ \nu(A\times A_1\times A_2).$$
Мы уже знаем в силу   ${ \bf Mix(1)}$, что для всех $A_1,A_2\in\B$ 
выполнено 
$$\nu(X\times A_1\times A_2)=\mu(A_1)\mu(A_2).\eqno (1)$$
Нам надо показать, что $\nu=\mu\times\mu\times\mu$.
Убедимся в том, что для $S\in H(T)$ имеет место равенство
$$ \nu(SA\times A_1\times A_2)=\nu(A\times A_1\times A_2). \eqno (2) $$
Действительно, 
$$\nu(SA\times A_1\times A_2) =\lim_m \mu(SA \cap T^{i_m}A_1\cap T^{j_m}A_2)= $$
$$ =\lim_m \mu(A \ \cap \ S^{-1}T^{i_m}A_1\ \cap \ S^{-1}T^{j_m}A_2) =$$
$$
=\lim_m \mu(A \cap T^{i_m}T^{-i_m}S^{-1}T^{i_m}A_1\cap T^{j_m}T^{-j_m}S^{-1}T^{j_m}A_2) =
$$
$$
=\lim_m \mu(A \cap T^{i_m}A_1\cap T^{j_m}A_2)= \nu(A\times A_1\times A_2).
$$
 Если  $\mu(A)=\mu(B)$, то множества $A$ и $B$  равносоставлены  посредством 
преобразований из  эргодической группы
$H(T)$. Поэтому с учетом доказанного свойства $(2)$  получаем
$$ \nu(A\times A_1\times A_2)=\nu(B\times A_1\times A_2).$$
Учитывая (1), окончательно приходим к 
$$\nu=\mu\times\mu\times\mu.$$
Тем же приемом  свойство ${ \bf Mix(k)}$ устанавливается для всех $k$. Это дало 
новое доказательство известных фактов о кратном перемешивани 
для гауссовских и пуассоновских динамических систем, см.  \cite{19}.

Замечание. Если действие некоторой группы  не обладает свойством кратного перемешивания, то оно, как мы показали,
 обладает бедной гомоклинической группой. Для действия Ледрапье на гармонической группе, как заметил С.В. Тихонов,  
 группа Гордина тривиальна: состоит из одного элемента $I$.


\section{ Мозаики и лемма Рохлина без $\eps$.}
В 1985 году А.М. Степин  предложил автору подумать над обобщением  результата  Лерера и Вейса \cite{LW}. Речь шла о следующем.  Пусть $T$ -- обратимое сохраняющее меру апериодическое преобразование вероятностного пространства
(например, отрезок $X$ единичной длины с мерой Лебега $\mu$). Лемма Рохлина-Халмоша утверждает, что \it
для любого $h$ и любого $\eps>0$  найдется разбиение  всего пространства на измеримые непересекающиеся множества
$$B, TB,T^2B,\dots, T^{h-1}B, \hat B,$$ причем  
$\mu(\hat B)<\eps.$  \rm Набор непересекающихся множеств вида  $B, TB,\dots, T^{h-1}B$ или  объединение этих множеств 
называют башней высоты $h$, а множество
$\hat B$ иногда называют  крышей над башней. 

Если потребовать эргодичность степени $T^n$ 
(отсутствие $T^n$-инвариантных нетривиьных множеств), то, как показали авторы \cite{LW}, 
 \it крышу над такой башней можно
поместить внутрь любого наперед заданного множества положительной меры. \rm

Есть ли аналог этой леммы для действия группы $\Z^2$, порожденного
 коммутирующими преобразования $R,S$ такими, что все $R^iS^j$ эргодичны при $i+j>0$? Уточним вопрос:
\it верно ли, что для любого $Y$, $\mu(Y)>0$, найдется башня  вида $U=\uu_{0\leq i,j <h}  R^iS^jB$
такая, что крыша $X\setminus U$ лежит в $Y$ ? \rm   

\bf Мозаика из квадратиков. \rm Оказалось, что ответ отрицательный \cite{88},\cite{92}, получить его помогает следующая  совсем детская игра:
\it
имеется   набор  красных квадратиков размера  $3\times 3$ и набор синих квдратиков размера 
 $1\times 1$, из них надо составить мозаику при   условии, что синие квадратики не касаются друг друга.
 \rm

\begin{center}
\includegraphics{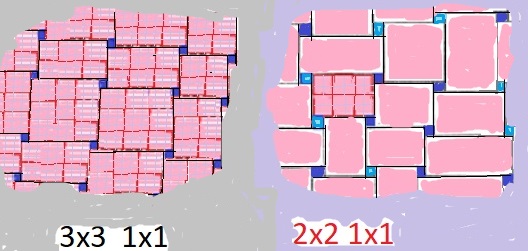} 
\end{center}

В полученном  мозаичном рисунке     наблюдается глобальная закономерность, хочется сказать  ''квазикристалличность''.
Энтропия   таких систем нулевая (мало  беспорядка), 
следовательно, для действий с положительной энтропией лемма Рохлина-Халмоша без  $\eps$ не выполняется.

 Мозаика из квадратиков $2\times 2$ и $1\times 1$  отличается  от   квазикристаллов,
она  аморфна.  Синие квадратики обладают    спинами $+1$ и $-1$, в то время как у квазикристаллов  
все синие квадраты имеют или спин $+1$, или все они  обладают спином $-1$.  Можно  сравнить  картинки  и  догадаться, что  такое спин маленьких квадратиков. 

\begin{center}
 \includegraphics{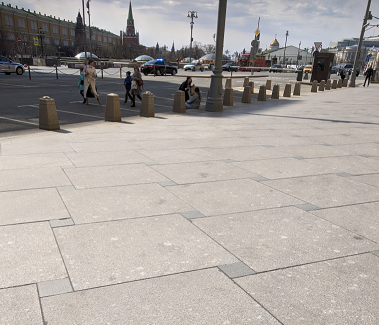}
\end{center}
\begin{center}
{\large \it  Первый тип  мозаик   используют в Москве, \\ например, на ул. Охотный Ряд.\rm}
\end{center}

\vspace{3mm}
\bf Нерешенная задача. \rm В прошлом тысячелетии  Б.Вейс прислал автору задачу а-ля  лемма Рохлина. 
 Сформулируем ее частный случай.  Рассмотрим    свободную группу $F_2$ с двумя образующими $a,b$ и ее стандартное бернуллиевское действие на пространстве $X=2^{F_2}$ с мерой Бернулли $\mu$ типа 
$(\frac 1 2, \frac 1 2 )$. 
Напомним, что   множество последовательностей из нулей и единиц, пронумерованных элементами группы $F_2$ и 
принимающих заданные  значения  на конечном  множестве $S\subset F_2$, по определению имеет меру  ${2^{-|S|}}$. 
Продолжая меру по Лебегу, получаем вероятностное пространство, 
его точки обозначаются через $x, \ x: F_2\to \{0,1\}$.
 Обратимые преобразования $T_a,T_b:X\to X$, определенные формулой
$$T_ax(g)=x(ag),\ g\in F_2, $$   сохраняют бернуллиевскую меру. 

\vspace{2mm}\it
Какова верхняя грань множества значений  $\mu(B)$ для таких измеримых  $B$, что 
множества $B,\ T_aB,\ T_bB,\ T_{a^{-1}}B, \ T_{b^{-1}}B $  не пересекаются?\rm

\vspace{2mm}
Наши знания о $\sup \mu(B)$ весьма скромны:   $\sup \mu(B) \in  [\frac {1} {\approx 17},\frac 1 5]$. 


Указанный набор  множеств образует  башню Рохлина  с конфигурацией 
$$\{a,b,e, a^{-1}, b^{-1}\}.$$

Заметим, что непересекающимися сдвигами этой конфигурации (крестами)   можно замостить всю группу $F_2$.

\begin{center}
\includegraphics{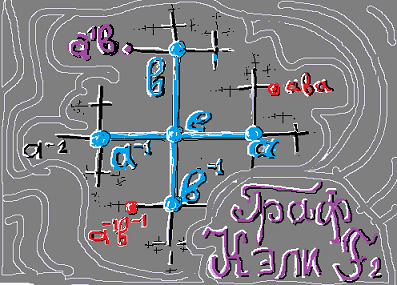}
\end{center}

 \section{Необычные и типичные свойства преобразований}
Говорят, что охраняющее меру обратимое преобразование
$ T: X \to X $ пространства Лебега $ (X, \mu) $ обладает {\it рангом один}, 
если существует последовательность 
$$ \xi_j= \{E_j, \ TE_j, \ T^2 E_j, \ \dots, T^{ h_j} E_j, \tilde {E}_j \} $$ 
измеримых разбиений пространства $ X $ таких, что 
$ \xi_j $ стремится к разбиению на точки.  Это означает, что любое множество 
конечной меры сколь угодно точно для всех достаточно больших $j$  приближается $ \xi_j$-измеримыми  множествами. 
Аналогично определяется ранг 1 для потоков (действий группы $\R$).

Действия ранга допускают конструктивное описание. Варьируя конструкцию, можно управлять ее свойствами, в частности,
спектром, т.е. наделять  спектральную меру действия  нужным  свойством.
Это направление, связанное с именами     А.Б. Катка, 
В.И. Оселедца, Д. Орнстейна, Д. Рудольфа, А.М. Степина, Р. Чакона и многих других, 
привело к обнаружению  новых эффектов в эргодической теории \cite{20}.
 
Рассмотрим свойство меры $\sigma$ на прямой $\R$,  которое неформально описывается так:   \it ''мера $\sigma\times\sigma$ прозрачна для некоторых  направлений, но   дает полную тень вдоль других.'' \rm  
 Подразумевается, что    произведение $\sigma\times\sigma$  под некоторыми углами проектируется однократно (биективно 
$mod \, 0$) в сингулярные меры   на некоторой  прямой  $L\subset \R^2$,  
но при этом есть  множество других углов, для которых  $\sigma\times\sigma$ проектируется  на $L$ с
бесконечной кратностью   в  меры, эквивалентные мере Лебега.
На роль  такой меры $\sigma$  выбирается  спектральная  мера специально построенного  потока  ранга один. 

\bf Парадоксальный эффект    Катка.  \it Плоскость $\R^2$  допускает
расслоение на  кривые, являющиеся  графиками непрерывных монотонных функций на $\R$, такое, 
что каждая кривая пересекалается  с  некоторым подмножеством плоскости $Y$
в точности по одной точке, причем  дополнение к  $Y$ имеет нулевую меру Лебега. \rm   (см. \cite{M}.) 

 Такое  расслоение можно изготовить, в частности, при помощи  
 упомянутой  нами меры $\sigma$ на $\R$. Надо  воспользоваться  тем, что произведение $\sigma\times\sigma$ 
проектируется однократно  в сингулярную меру на прямой. 

\bf Контролируемое отклонение от перемешивания.  \rm  В эргодической теории хорошо известен факт о том, что если преобразование  перемешивает  на  некотором (редком) подмножестве моментов времени, то оно необходимо перемешивает на некотором (густом) множестве  плотности 1.  

Задача  Виталия Бергельсона:  \it пусть для  измеримых подмножеств $A,B$ вероятностного 
пространства выполняется $$\mu(T^{n^2}A\cap B) \to \mu(A)\mu(B), n\to\infty, $$ верно ли, что 
$$\mu(T^{n}A\cap B) \to \mu(A)\mu(B), n\to\infty \ ?$$  \rm (Этот вопрос недавно сообщил автору Эль Абдалаоуи.)

Отрицательный ответ дает следующее общее утверждение.

\vspace{3mm}
\bf Теорема. \it Пусть $I_j=[a(j), b(j)]$ -- последовательность отрезков, причем  $b(j)-a(j)\to \infty$ и  $a(j+1)>b(j)$ для всех $j\in\N$.
Найдется   неперемешивающее преобразование  такое, что всякая его неперемешивающая последовательность содержит подпоследовательность, лежащую  в объединении отрезков $I_j$.

\vspace{2mm}
\bf Следствие. \it Для любого бесконечного множества $M\subset\N$ нулевой плотности некоторое  неперемешивающее преобразование  перемешивает вдоль $M$.\rm

В работе \cite{19a}  предложен класс  конструкций ранга один (специальные перекладывания 
бесконечного семейства отрезков), которые заданы    натуральными параметрами $s_j$.
Некоторое представление о том, как они  устроены, дает  следующий  рисунок.

\begin{center}
\includegraphics{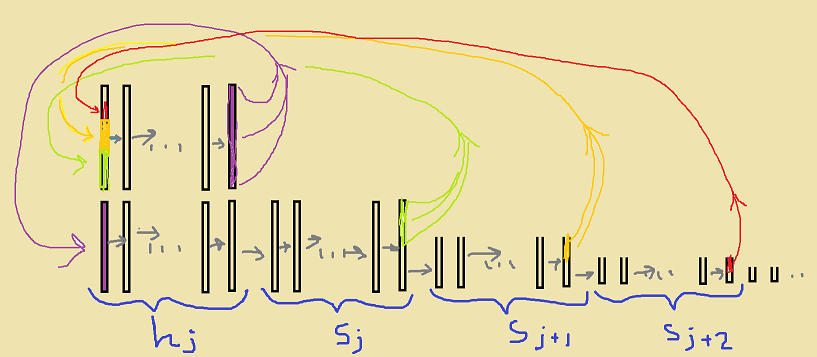}
\end{center}

Определим рекуррентно числа $h_j$ соотношением 
$h_{j+1}=2h_j+ s_j$ с условием  $s_j>>h_j$. Для наших конструкций выполняется отклонение от перемешивания вида
$$ \mu(T^{h_j}A\cap A) \to \mu(A)/2, j\to\infty, $$ 
где $A$ -- множество конечной меры. Но главная особенность  в том, что
мы знаем описание всех  неперемешивающих последовательностей.
Действительно, если для некоторых множеств $A,B$ конечной меры и последовательности $\{n(k)\}$ выполнено
$$ \liminf_{n(k)\to\infty} \mu(T^{n(k)}A\cap B) >c>0, $$
то, как показано в \cite{19a},  для всех больших $k$ ограничена последовательность 
$$ n(k)- h_{j_1(k)}\pm h_{j_2(k)}\dots \pm h_{j_m(k)},$$
где $m$ также  ограничено ($2^{-m}\mu(A)\geq c$) и $j_1>j_2>\dots>j_m$.

Как подобрать параметры конструкции для заданной последовательности отрезков $\{I_j\}$? Выберем   редкую
подпоследовательность отрезков, обозначая ее снова через $\{I_j\}$, и параметры конструкции $s_j$ так,
чтобы $h_j$ попадало в середину отрезка $I_j$, причем  длина отрезка $I_j$ была много больше $h_{j-1}>>h_{j-2}>>\dots$.
Тогда для упомянутой неперемешивающей последовательности $\{n(k)\}$  верно, что для всех больших $k$ будет выполняться
$n(k)\in I_{j_1(k)}$.

Мы доказали теорему для случая бесконечной меры. Пуассоновские или  гауссовские надстройки над этими преобразованиями
дают нужные примеры для случая вероятностной меры, так как они сохраняют свойства последовательности быть или не быть 
перемешивающей. Можно строить примеры непосредственно в пространстве конечной  меры, но оно тесное, а это   усложняет конструкцию.

\vspace{3mm}
\bf Проблема Банаха (\cite{Ulam}): 
\it существует  ли  сохраняющее меру Лебега  преобразование $T$  прямой $\R$ такое, 
что для некоторой  функции $f\in L_2(\R)$  множество
$\{f\circ T^n: n\in\Z\}$ является  полной ортонормированной системой в $L_2(\R)$? \rm 

Конструкции ранга один с бесконечной инвариантной мерой  -- естественные кандидаты на роль таких преобразований.
 
 \bf Типичные свойства.   \rm Фиксируя  вероятностное пространство $(X,\B,\mu)$ и  плотное в  $\B$ семейство
 множеств  $\{A_i\}$, снабдим 
группу  $ Aut=Aut(X,\B,\mu)$  полной метрикой Халмоша $\rho $: 
$$ \rho(S,T)=\sum_i 2^{-i}\left(\mu(SA_i\Delta TA_i)+\mu(S^{-1}A_i\Delta T^{-1}A_i)\right).$$
 Говорят, что свойство преобразований типично, если множество всех элементов группы $Aut$,
  обладающих этим свойством,   содержит  счетное пересечение открытых множеств, 
плотных в  $Aut$.

\bf Некоторые свойства типичного преобразования. \it
 Преобразование не сопряжено своему обратному, не перемешивает, имеет нулевую энтропию,
является слабо перемешивающим,  обладает рангом один и   не допускает 
 нетривиальных самоприсоединений с парной независимостью.

Преобразование включается в поток, спектральная мера $\sigma$ которого такова, что  $\sigma\times\sigma$ 
 проектируется однократно   на типичную  прямую   в $\R\times\R$.

\bf Нерешенные задачи о типичных свойствах. \rm Современное состояние теории типичных преобразований отражено в
работе \cite{GTW}.  Мы ограничимся упоминанием двух   вопросов, возникших у автора в связи изучением гомоклинических групп \cite{19} и исследованием типичных свойств   энтропии  относительно последовательности конечных конфигураций 
в группе \cite{21}.
 
Эргодические преобразования $T$ с положительной энтропией  обладают   свойством: \it для некоторого
 апериодического преобразования $S$ найдется  последовательность $n(i)\to\infty$ такая, что выполнено
$$T^{-n(i)}ST^{n(i)}\to I.$$
\it Является  ли это свойство типичным?\rm

Обозначим     $$K^{Aut}=\{J^{-1}SJ:  S\in K, J\in Aut\}.$$ 
 \it Верно ли, что  для всякого  компакта $K\subset Aut$  множество $Aut\setminus K^{Aut}$  типично?
\rm Ответ получен в случае, когда компакт $K$ лежит в классе преобразований с нулевой классической
энтропией \cite{21}.

 \large

\hfill vryzh@mail.ru
\end{document}